\documentclass[10pt,leqno]{article}
\usepackage{amsmath,txfonts}
\usepackage{graphics,hyperref}
\usepackage{showkeys}
\usepackage{color}
\numberwithin{equation}{section}

\newcommand{\eps}{\varepsilon}
\newcommand{\dt}{\delta}
\newcommand{\om}{\omega}

\newcommand{\la}{\lambda}
\newcommand{\bt}{\beta}
\newcommand{\ga}{\gamma}

\renewcommand{\L}{\varmathbb{L}^{2}}
\newcommand{\E}{\mathcal{E}}
\newcommand{\I}{\mathcal{I}}
\renewcommand{\H}{\varmathbb{H}^{1}}
\newcommand{\M}{\mathcal{M}_{\ga}}
\newcommand{\Mdue}{\mathcal{M}_{(\delta_{1},\delta_{2})}}
\newcommand{\Md}{\mathcal{M}_{\delta}}
\newcommand{\Mdd}{\mathcal{M}_{(\delta,\delta)}}
\newcommand{\G}{\mathcal{G}}

\newcommand{\Se}{\mathcal{S}}
\newcommand{\B}{\mathcal{B}}
\newcommand{\A}{\mathcal{A}}
\renewcommand{\P}{\mathcal{P}}
\newcommand{\cdue}{c_{\gamma}}

\newcommand{\mnehari}{m_{\nh}}
\newcommand{\nh}{\mathcal N}
\newcommand{\zbo}{z_{\beta}^{\omega}}

\newcommand\R{\varmathbb{R}}
\newcommand\C{\varmathbb{C}}
\newcommand\huno{H^{1}}
\newcommand\ldue{L^{2}}

\newcommand{\rn}{\R^{n}}

\newcommand{\pa}{\partial}
\newcommand{\lapl}{\Delta}

\newtheorem{theorem}{Theorem}[section]
\newtheorem{remark}[theorem]{Remark}
\newtheorem{lemma}[theorem]{Lemma}
\newtheorem{proposition}[theorem]{Proposition}
\newtheorem{corollary}[theorem]{Corollary}
\newtheorem{definition}[theorem]{Definition}
\newcommand{\bos}{\begin{remark}\rm}
\newcommand{\eos}{\end{remark}}
\newcommand{\bte}{\begin{theorem}}
\newcommand{\ete}{\end{theorem}}
\newcommand{\bpr}{\begin{proposition}}
\newcommand{\epr}{\end{proposition}}
\newcommand{\bdf}{\begin{definition}\rm}
\newcommand{\edf}{\end{definition}}
\newcommand{\bco}{\begin{corollary}}
\newcommand{\eco}{\end{corollary}}
\newcommand{\ble}{\begin{lemma}}
\newcommand{\ele}{\end{lemma}}
\newcommand{\bdm}{\begin{displaymath}}
\newcommand{\edm}{\end{displaymath}}
\newcommand{\beq}{\begin{equation}}
\newcommand{\eeq}{\end{equation}}
\newcommand{\bdim}{{\noindent{\bf Proof.}\quad}}
\newcommand{\edim}{{\unskip\nobreak\hfil\penalty50
 \hskip2em\hbox{}\nobreak\hfil\mbox{\rule{1ex}{1ex} \qquad}
  \parfillskip=0pt \finalhyphendemerits=0\par\medskip}}

\newcommand{\ov}{\overline}

\newcommand{\dyle}{\displaystyle}
\newcommand{\dys}{\displaystyle }

\begin{document}

\title{Orbital stability property for\\ 
coupled nonlinear Schr\"odinger equations}
\author{Liliane Maia\thanks{Research partially 
supported by Projeto Universal/CNP$_{q}$ and FAPDF.},\; 
Eugenio Montefusco\thanks{Research 
supported by MIUR project {\it Metodi Variazionali ed 
Equazioni Differenziali non lineari}.},\; Benedetta Pellacci$^{\dag}$}
\date{}
\maketitle

\begin{abstract}
Orbital stability property for weakly coupled nonlinear Schr\"odinger equations 
is investigated. Different families of orbitally stable standing waves solutions 
will be found, generated by different classes of solutions
of the associated elliptic problem. In particular,
orbitally stable standing waves can be generated by least action solutions, but
also by solutions with one trivial component whether or not they are ground states. Moreover,  standing waves  with components 
propagating with the same frequencies are orbitally stable if
generated  by vector  solutions of a  suitable single Schr\"odinger weakly coupled system, even if they are not ground states.
\end{abstract}

\section{Introduction}\label{sec.introd}
We consider the following Cauchy problem for two
coupled nonlinear Schr\"odinger equations
\beq\label{schr}
\begin{cases}
i\pa_{t}\phi_{1}+\lapl\phi_{1}+\left(|\phi_{1}|^{2p-2}
  +\beta|\phi_{2}|^{p}|\phi_{1}|^{p-2}\right)\phi_{1}=0, & \\
i\pa_{t}\phi_{2}+\lapl\phi_{2}+\left(|\phi_{2}|^{2p-2}
  +\beta|\phi_{1}|^{p}|\phi_{2}|^{p-2}\right)\phi_{2}=0, & \\
\phi_{1}(0,x)=\phi_{1}^{0}(x), \quad \phi_{2}(0,x)=\phi_{2}^{0}(x),
\end{cases}
\eeq
where $\Phi=(\phi_{1},\phi_{2})$ and $\phi_{i}:\R\times\R^n\to\C$, 
$\phi_{i}^{0}:\R^n\to\C$, $p>1$ and $\beta$ is a real 
positive constant.
\\
Coupled nonlinear Schr\"odinger equations appear in the study 
of many physical processes. For instance, such equations with 
cubic nonlinearity model the nonlinear interaction of two wave 
packets, optical pulse propagation in birefringent fibers 
or wavelength-division-multiplexed optical systems 
(see  \cite{man,men},
\cite{aa,ka} and the references therein). 
\\
A soliton or standing wave solution is a solution of the form 
$\Phi(x,t)=(u_{1}(x)e^{i\om_{1}t},u_{2}(x)e^{i\om_{2}t})$ where 
$U(x)=(u_{1}(x),u_{2}(x)):\R^{n}\to \C^{2}$ 
is a solution of the elliptic system
\beq\label{ellittico}
\begin{cases}
  -\lapl u_{1}+\om_{1}u_{1}=\left(|u_{1}|^{2p-2}
  +\bt|u_{1}|^{p-2}|u_{2}|^{p}\right)u_{1}, & \\
  -\lapl u_{2}+\om_{2}u_{2}=\left(|u_{2}|^{2p-2}
  +\bt|u_{2}|^{p-2}|u_{1}|^{p}\right)u_{2}. &
\end{cases}
\eeq
Among all the standing waves we can distinguish between ground 
and bound states. A ground state corresponds to a least 
action solution $U$ of \eqref{ellittico}; while all the other critical 
points of the action functional give rise to bound states 
(or excited states) of \eqref{schr}. A ground state generates
a one-hump soliton of \eqref{schr} because it is 
nonnegative, radially symmetric and decades exponentially 
at infinity (\cite{bs}). On the other hand, 
vector multi-hump solitons are of much interest in the 
applications, for example they have been observed in 
photorefractive crystals \cite{msc}. 
\\
When investigating 
stability properties of a given set of solution, it is natural to take into 
account the rotation invariance of the problem, and this is done by the 
orbital stability. Roughly speaking, this means
that if an initial datum $\Phi^{0}$ is close to a ground state 
$U$ then all the orbit generated by $\Phi^{0}$ remains close to 
the soliton generated by $U$ up to translations or phase 
rotations.
\\
For the single Schr\"odinger equation it has been proved that
the orbital stability property is enjoyed by standing waves 
raised by least action solutions. This result can be deduced from
the two following facts (see \cite{cl} and Section 8 in \cite{caz}):
\begin{itemize}
\item[(a)]
every least action solution can be associated, by a bijective correspondence,
to a minimum point of the energy constrained to the $L^{2}$ sphere
with a suitable choice of the radius.
\item[(b)]
Conservation laws and compactness properties of this minimization
problem imply that the set of minimum points of the energy on this sphere 
manifold generates stable standing waves.
\end{itemize}

Moreover, in \cite{gss,ss} it is proved that every  critical 
point of the action with Morse index larger than one give rise to instability.
Taking into consideration the result of \cite{k} the stability of 
the standing wave $e^{i\om t}z_{\om}$ of the single Schr\"odinger equation 
holds if and only if $z_{\om}$ is the minimum point of the associated 
energy functional constrained to the $L^{2}$ sphere of radius 
$\|z_{\om}\|^{2}_{L^{2}}$.

This and (a) are the reasons why ground states are the most desirable solution 
for the single Schr\"odinger equation. 

With this situation in mind, large effort has been done in the last 
few years to find ground states of \eqref{ellittico}. 
In \cite{ka,shm} numerical arguments or analytical 
expansions have been employed to produce different 
families of solitons. 
The investigation has been improved by means of 
variational methods. In \cite{ac,bw,dfl,lw,mmp,si})
assumptions on the constant $\bt$ are stated in order 
to distinguish between ground states with both nontrivial 
components ({\it vector ground states}) and ground states 
with one trivial component ({\it scalar ground states}). 
It has been discovered that there exists vector ground states for the constant $\bt$ sufficiently large in dependence on the frequencies ratio, while if $\bt$ is small least action solutions have necessarily one trivial component. Moreover, in \cite{ac} it is clarified the 
difference between scalar and vector positive solutions in dependence to 
different geometrical properties of the action functional. For
$\bt$ small the scalar ground states are critical points of the action functional with Morse index equal to one, while for $\bt$ large 
these kind of solutions have larger Morse index. 
Since stable standing waves should be generated only by
ground states, these results suggested the idea that stable standing waves
should be given by scalar solutions  for $\bt$ small and by vector ground
states for $\bt$ large.
This opinion is confirmed also in \cite{l} where this topic has 
been studied  for different evolution systems, and the orbital
stability property is shown to be enjoyed by standing 
waves associated to solution of the corresponding
elliptic system with Morse index equal to one.
\\
For the cubic NLS systems in the one dimensional case,
this subject has been recently studied in some interesting papers using 
numerical and analytical methods. 
In \cite{y} it is conjectured, based on numerical 
evidence, that single-hump soliton are stable while
multi-hump vector solitons are all linearly 
unstable and this is proved by numerical and analytical 
arguments in \cite{py} for $p=2$ in \eqref{schr} and 
for any $p$ for special families of multi-hump vector solitons. 
In \cite{pk} a stability criterion is found to study the 
stability property of some families of single-hump vector 
solitons. 
When tackling this matter for weakly coupled Schr\"odinger 
equations by means of variational methods, one has to take
into account that the $L^{2}$ norms of the components 
are conserved se\-pa\-rately (see \cite{fm,o}). 
So that we can consider different constrains on which minimize the energy. When we choose the sphere with respect of
the  $L^{2}\times L^{2}$ norm, we obtain ground states, however we do not know  whether or not they are scalar or vector solutions. Otherwise, we could
try to mi\-ni\-mize the energy constraining the $L^{2}$ norms separately,
this approach will permit us to know in advance if we will find scalar or vector
solutions even if they may be not least action solutions. 
The first approach consists in solving the minimization problem
\beq\label{intromini}
\E(u)=\inf_{\M}\E\quad \text{where}\quad\M=\left\{ U=(u_{1},u_{2})\in\huno\times\huno :
\om_{1} \|u_{1}\|_{\ldue}^{2}+
\om_{2} \|u_{2}\|_{\ldue}^{2}=\ga\right\}
\eeq
and 
$$
\E(U)=\E(u_{1},u_{2})=\frac12\|\nabla u_{1}\|_{2}^{2}+
\frac12\|\nabla u_{2}\|_{2}^{2}-\frac{1}{2p}
\left(\|u_{1}\|_{2p}^{2p}+\|u_{2}\|_{2p}^{2p}
+2\beta\|u_1 u_2\|_{p}^{p}\right),
$$
for $1<p<1+2/n$ in order to have global existence of \eqref{schr}
(see \cite{fm}).\\
For a suitable choice of $\gamma$, we will find that this problem 
has a solution corresponding, in  a bijective correspondence, 
to ground states of  \eqref{ellittico}. Moreover, the set of the 
solutions generates stable sanding waves (see Theorem \ref{mainstab}). 
This conclusion is in accordance to the single equation case, since 
these solutions have Morse index equal to one. \\
When adopting the second approach, we are naturally lead to the minimization problem
\beq\label{intromini2}
\E(u)=\!\!\!\inf_{\Mdue}\!\!\!\E\quad\text{where}\quad\Mdue=\left\{ U=(u_{1},u_{2})\in\huno\times\huno : \|u_{1}\|_{\ldue}^{2}=\delta_{1},\,\|u_{2}\|_{\ldue}^{2}=\delta_{2}\right\}.
\eeq
When  $\delta_{2}$ (or $\delta_{1}$) is equal to zero we obtain as 
minimum point the couple $(z_{\om_{1}},0)$ (or $(0,z_{\om_{2}})$) 
where $z_{\om_{1}}$ ($z_{\om_{2}}$) is the unique  positive 
solution of the first (second) equation in \eqref{ellittico}, recall that this solution is a ground state of \eqref{ellittico} only for $\beta$ small. 
However, we will show that they still produce orbitally stable standing 
waves for any $\beta>0$ (see Theorem \ref{mainstab2}).
This result is in accordance with the conjecture in \cite{y}.
But they are in contrast with the expectation that only ground states should give rise to orbitally stable waves, since they have Morse index greater than one for $\beta$ large. 
The case $\delta_{1}=\delta_{2}$ has been tackled in \cite{oh}
for $p=2,\,n=1$ and $\beta=1$, and it is proved that the set of solutions 
of \eqref{intromini2} give rise to orbitally stable solutions of \eqref{schr}.
Here we will extend the result in \cite{oh} for higher dimension and for 
every $\beta>0$ see Theorem \ref{mainstab3}. Moreover, as in \cite{oh}, we will show that, for a suitable choice of $\delta$ (and $\delta_{1}=\delta_{2}=\delta$) the set of solutions of \eqref{intromini2} is given by
$$
\B=\left\{(e^{i\theta_{1}}\zbo(\cdot-y),e^{i\theta_{2}}
\zbo(\cdot-y)),\,\theta_{1},\,\theta_{2}\in \R,\, y\in\R^n\,\right\},
$$
where $\zbo$ is the unique positive solution of the problem
$$
\begin{cases}
-\Delta u+\om u=(1+\beta)|u|^{2p-2}u &\text{in $\R^{n}$},
\\
\;\; \;u(x)\to 0 &\text{as $|x|\to \infty$},
\end{cases}
$$
Moreover, arguing as in \cite{oh} we will demonstrate that $\B$
also characterizes the set of least action  solution 
of \eqref{ellittico} when  $\omega_{1}=\omega_{2}$ and when one 
prescribes both of the components to be different from zero
(see Theorem \ref{carabound}). So that our result provide a complete characterization of the set of solutions found in Theorem 1 in \cite{si}
and in Theorem 2 in \cite{lw} for $\lambda_{j}=\om$ for every $j$.
Let us stress again that the set $\B$ is made of ground states only for
$\beta\geq 1$.
Then, for $\bt $ large we have at least two families of  orbitally stable solution of \eqref{schr}, the ones generated by ground states, which we know have
both nontrivial components and the ones produced by the scalar solution.
We can reach the same conclusion for any $\beta$ but for $\omega_{1}=\omega_{2}$: orbital stability is enjoyed by the standing waves generated 
by scalar solutions and by vector solution solutions of \eqref{intromini2}, 
the former are ground states for $\bt$ small, the latter are
ground states for $\bt$ large. Unfortunately we cannot handle the case 
$\delta_{1}\neq\delta_{2}$, and to our knowledge the question of whether or not the set of solution of \eqref{intromini2} gives rise of orbitally stable solutions for any $\delta_{1},\,\delta_{2}$ is open.
\\
Finally, adapting the arguments in \cite{caz},
we will also show an instability result in the supercritical 
case $p=1+2/n$, for ground state solutions, scalar solutions 
and for the set $\B$ (for $\om_{1}=\om_{2}=\om$),
as a consequence of blowing up in finite time.
While, for the critical case the instability is produced by every solution
of \eqref{ellittico}.
\\
The paper is organized as follows. In Section \ref{setting} 
we state our main results. The definitions and preliminary 
results, preparatory to the proofs, are presented in section \ref{preliminary}. 
In section \ref{proofs} we give the proofs of our main results.
A section of conclusion  comments the results obtained.

\section{Setting of the problem and main results}\label{setting}

Our analysis will be carried out in the functional spaces
$\L=\ldue(\rn,\C)\times\ldue(\rn, \C)$ and
$\H=\huno(\rn,\C)\times\huno(\rn,\C)$. 
We recall that the inner product between $u,\,v\in \C$ is given by 
$u\cdot v=\Re(u\ov{v})=1/2(u\ov{v}+v\ov{u})$. Then for
$\om=(\om_{1},\om_{2})$, $\om_{i}\in \R,\,\om_{i}>0$, we can define
an equivalent inner product in $\L$ given by
\bdm
(\Phi|\Psi)_{\om}=\Re\int\left[\om_{1}\phi_{1}\ov{\psi}_{1}
+\om_{2}\phi_{2}\ov{\psi}_{2}\right],\qquad \forall\, 
\Phi=(\phi_{1},\phi_{2}),\, \Psi=(\psi_{1},\psi_{2})
\edm
and an equivalent norm
\bdm
\|\Phi\|_{2,\om}^{2}= \|\phi_1\|^2_{2,\om_{1}}+\|\phi_2\|^2_{2,\om_{2}},
\quad \text{where}\quad
\|\phi_i\|^2_{2,\om_{i}}=\om_i\|\phi_i\|^2_{2}
=\om_{i}\int\phi_{i}\ov{\phi}_{i}
\edm
for $i=1,2$. It is known (see  Remark 4.2.13 in \cite{caz}) that \eqref{schr} is 
locally well posed in time, for $p<n/(n-2)$ when $n>2$ and for 
any $p$ for $n=1,2$, in the space $\H$ endowed with the norm 
$\|\Phi\|_{\H}^{2}=\|\nabla\Phi\|_{2}^{2}+\|\Phi\|_{2,\om}^{2}$ 
for every $\Phi=(\phi_1,\phi_2)\in\H$. Moreover we set
the ${\varmathbb{L}^{p}}$ norm as
$\|\Phi\|_p^p=\|\phi_1\|_p^p+\|\phi_2\|_p^p,$ 
for $p\in\left[1,+\infty\right)$. It is well known that the masses
of the components of a solution and its total energy are preserved in time, 
that is the following conservation laws hold  (see \cite{fm,o}):
\beq\label{mass}
\|\phi_1\|_2^2=\|\phi^{0}_{1}\|_2^2,\qquad
\|\phi_2\|_2^2=\|\phi^{0}_{2}\|_2^2,
\eeq
\beq\label{energy}
\E\left(\Phi(t)\right)=\frac{1}{2}\left\|\nabla\Phi (t)\right\|_2^2 
-F\left(\Phi(t)\right) =\frac{1}{2}\left\|\nabla\Phi^{0}\right\|_2^2
-F\left(\Phi^{0}\right)=\E(0),
\eeq
where
\beq\label{defF}
F(\Phi)=\frac{1}{2p}\left(\|\Phi\|_{2p}^{2p}
+2\beta\|\phi_1\phi_2\|_{p}^{p}\right).
\eeq
In \cite{fm} it is proved that the solution of this Cauchy 
problem exists globally in time, under the assumption
\beq\label{pzero}
p<1+\frac2n.
\eeq
In order to study orbital stability properties, 
we will use the functional energy (see \eqref{energy}) and the action 
\bdm
\I(U)=\dfrac{1}{2}\|U\|_{\H}^{2}-F(U)=\E(U)+\frac12\|U\|_{2,\omega}^{2}.
\edm


\bdf\label{ground}
We will say that a ground state solution $U$ of 
\eqref{ellittico} is a solution of the following 
minimization problem
\begin{equation}\label{cne}
\I(U)=m_{\nh} :=\inf_{\nh}\I(W) \quad \text{where}\quad
\nh\,:\,=\{W\in \H\setminus\{0\}\,:\,\langle \I'(W),W\rangle =0\,\}.
\end{equation}
$\nh$ is called in the literature Nehari manifold 
(see \cite{mmp,si}). Moreover, we will denote with 
$\G$ the set of the ground state solutions.
\edf

\bdf\label{bound}
We will say that a positive bound state solution $U$ 
of \eqref{ellittico} is a solution of the following 
minimization problem
\begin{equation}\label{cnedue}
\begin{array}{l}
\dyle \I(U)=m_2 :=\inf_{\nh_{2}}\I(W) 
\quad\hbox{ where} \\ 
\nh_{2}\,:\,=\{W=(w_{1},w_{2})\in \H \,:\, w_{1},w_{2}\neq0,
\, \langle \partial_{1}\I(W),w_{1}\rangle =
\langle \partial_{2}\I(W),w_{2}\rangle = 0\,\},
\end{array}
\end{equation}
where $\partial_{1}\I(W)$ ($\partial_{2}\I(W)$) is the partial 
derivative with respect to the first (second) component.
$\nh_{2}$ is called in the literature Nehari set 
(see \cite{lw,si}). Moreover, we will denote with 
$\B$ the set of such bound state solutions.
\edf
It is well known (see Section 8 in\cite{caz}, \cite{k}) that all
the solutions of  the elliptic problem
\beq\label{schrequa}
\begin{cases}
-\Delta u+\om u=(1+\beta)|u|^{2p-2}u &\text{in $\R^{n}$},
\\
\;\; \;u(x)\to 0 &\text{as $|x|\to \infty$},
\end{cases}
\eeq
for $\om>0$ and $\bt\geq0$, are given by $v(x)=e^{i\theta}\zbo(x-y)$
where $\theta\in \R$ $y\in \R^{n}$ and $\zbo$ is the  unique 
least energy solution , where $\zbo$ is the unique positive
least energy solution in $H^{1}(\R^n,\R)$ of \eqref{schrequa}.
Let us recall that
\beq\label{defzbo}
\zbo(x)= \left(\dfrac{\om}{1+\bt}\right)^{1/2(p-1)} 
\!\!\!\! z^{1}_{0}\left(\sqrt{\om}x\right).
\eeq
\bdf\label{scalari}
Problem \eqref{ellittico} admits also scalar solutions,
$U=(u_{1},0)$ (or $(0,u_{2})$).  We will denote with 
$\Se$ the set of such solutions.
The uniqueness result in \cite{k} for the single Schr\"odinger equation, 
gives us the following characterization for the set $\Se$.
\bdm
\Se =\left\{(e^{i\theta}z^{\om_{1}}_{0}(\cdot-y),0),\;\theta\in \R,\,
y\in\R^n\right\}\cup
\left\{(0,e^{i\theta}z^{\om_{2}}_{0}(\cdot-y)),\;\theta\in \R,\,
y\in\R^n\right\}
\edm
where $z^{\om}_{0}$ is defined in \eqref{defzbo}
(with $\om_{1}=\om$ or $\om_{2}=\om$).
\edf
\bos
The results contained in \cite{ac,bw,dfl,mmp,lw,si}
show that, depending on the parameters $\om_{1},\om_{2},\bt$, 
the set $\G$ may coincide with either $\B$ or $\Se$. 
\\
For $\bt$ sufficiently large in dependence on $\om_{1}, \om_{2},$
$\G=\B$  and the point in $\Se$ are scalar bound states solutions.
While, $\G=\Se$ for $\bt$ small.\\
In the particular case  $\om_{1}=\om_{2}=1$ ground states have 
both nontrivial components if and only if $\beta\geq 1$ (see \cite{ac,bw,mmp,si}), 
so that for $\beta\geq1$ $\G=\B$, while for $\beta< 1$ $\G=\Se$.
\eos

Let us recall the orbital stability property for a set of solutions 
$\mathcal{F}$, introduced for the single equation case in \cite{cl}.

\bdf 
A set $\mathcal{F}\subseteq \H$ of solutions of Problem \eqref{ellittico} 
is orbitally stable if for any 
$\eps>0$  there exists $\dt=\dt(\eps)>0$ 
such that
\bdm
\text{if } \inf_{U\in\mathcal{F}} \| \Psi^{0}-U\|_{\H}< \dt,
\quad \text{ then }\quad
\sup_{t\geq0} \inf_{V\in\mathcal{F}} \|\Psi(t,\cdot)-V\|_{\H}<\eps,
\edm
where $\Psi$ is the global solution of \eqref{schr} with initial
datum $\Psi^{0}$.
\edf
\bos
We call the property in the previous definition orbital stability of
${\mathcal F}$ because every element $(u_{1},u_{2})$ of ${\mathcal F}$ 
generates an orbit given by the standing wave 
$(e^{i\om_{1}t}u_{1},e^{i\om_{2}t}u_{2})$.
\eos
Roughly speaking, a set $\mathcal{F}$ is orbitally stable 
if any orbit generated from an initial datum $\Psi^{0}$ 
close to an element of $\mathcal{F}$ remains close to 
$\mathcal{F}$ uniformly with respect to the time.

Up to now the uniqueness of the ground state solution is an open
problem for system \eqref{schr}, so that in the definition of
orbital stability we have to take into account the possibility
of a solution $\Psi$ to go from a ground state $U$
to a different ground state solution $V$; with this respect, 
it would be very interesting to know, at least, if ground states 
are isolated. In addition, we will show that there exist also 
other sets of orbitally stable solutions, then our definition 
has to take into account this aspect. 

Our main results are the following ones.

\bte\label{mainstab}
Assume \eqref{pzero}. For any $\bt,\om_1,\om_2>0$ the set
$\G$ is orbitally stable.
\ete

Different from the single equations case, we have other 
families of orbitally stable solutions for the system, 
as the next results show.

\bte\label{mainstab2}
Assume \eqref{pzero}. For any $\bt,\om_1,\om_2>0$ the set
$\Se$ is orbitally stable.
\ete

\bos
The preceding results imply that the problem \eqref{schr}
possesses at least two family of orbitally stable standing 
waves for $\bt$ sufficiently large, ground state standing 
waves and scalar ones. While, for $\bt$ small the stability 
property of scalar standing waves is a consequence both of 
Theorems \ref{mainstab} and \ref{mainstab2}, since $\G=\Se$ 
in this case.
\eos

If $\om_{1}=\om_{2}=\om$ the set $\B$ is completely characterized 
in the next result.

\bte\label{carabound}
Assume $\om_{1}=\om_{2}=\om>0$.
For any $\bt\geq0$ it holds
$$
\B =\left\{(e^{i\theta_{1}}\zbo(\cdot-y),e^{i\theta_{2}}
\zbo(\cdot-y)),\,\theta_{1},\,\theta_{2}\in \R,\, y\in\R^n,\,\zbo 
\text{defined in \eqref{defzbo}}\right\}
$$
In other words the set $\B$ is described by the standing wave 
of the single equation, up to translations and phase shifts
of the components.
\ete

This characterization of the set $\B$ leads us to show that 
the set $\B$ is orbitally stable even for $\beta$ small.

\bte\label{mainstab3}
Assume \eqref{pzero}. For any $\bt\geq0$ and $\om_1=\om_2=\om>0$ 
the set $\B$ is orbitally stable.
\ete

\bos
\begin{enumerate}
\item
These results imply that solutions that starts from initial data close 
to ground states with both nontrivial components  remain close to 
orbits generated by ground states with both nontrivial components. 
While, solutions that start close to $\Se$ will stay  close to orbits 
generated by $\Se$.
\item
>From the preceding results we deduce that, for $\om_{1}=\om_{2}$,
$\B$ and $\Se$ are  always orbitally stable sets independently of $\bt$. 
When $\om_{1}\neq \om_{2}$, we have that $\Se$ is always orbitally 
stable, while we can prove that $\B$ is orbitally stable only  
when it coincides with $\G$, that is for $\bt$ large. 
It is an open problem to study the stability property for 
$\B$ for any $\om_{1}\neq \om_{2}$. Our results  cover 
the following cases:
\begin{enumerate}
\item
$\om_{1}=\om_{2}$ positive; for any $\bt\geq 0$.
\item
$\om_{1}\neq\om_{2}$, $\bt$ large (in dependence of $\om_{1}/\om_{2}$)  
such that $\B=\G$.
\end{enumerate}
\end{enumerate}
\eos

We will also prove an instability result for the sets $\G,\, \Se$ and $\B$
(for $\om_{1}=\om_{2}$) in dependence of the exponent $p$. 
More precisely, we will show the following results.

\begin{theorem}\label{instap}
Assume $p<n/(n-2)$. For any $\om_{1},\om_{2},\bt>0$ 
the following conclusions hold:

a) Let $p>1+2/n$, then the sets $\G,\,\Se$ are unstable in the following sense: 

For any $U\in\G$ (or $U\in \Se$) and $\eps>0$ there exists 
$U_{\eps}^{0}$ with 
$\|U_{\eps}^{0}-U\|_{\H}\leq \eps$ such that the solution 
$\Phi_{\eps}$ satisfying $\Phi_{\eps}(0)=U_{\eps}^{0}$ 
blows up in a finite time in $\H$.

b) Let $p=1+2/n$, then every solution of \eqref{ellittico} 
is unstable in the sense of the previous conclusion.
\end{theorem}

\begin{theorem}\label{instapB}
Assume $1+2/n<p<n/(n-2)$ and $\om_{1}=\om_{2}=\om$. 
For any $\bt>0$, the set $\B$ is unstable in the following sense: 

For any $U\in\B$ and $\eps>0$ there exists 
$U_{\eps}^{0}$ with $\|U_{\eps}^{0}-U\|_{\H}\leq \eps$ 
such that the solution $\Phi_{\eps}$ satisfying $\Phi_{\eps}(0)=U_{\eps}^{0}$ blows up in a finite time in $\H$.
\end{theorem}
\section{Minimization Problems}\label{preliminary}

\subsection{Ground states}

In this section we will present some general results which 
will be useful in proving Theorems \ref{mainstab}, 
\ref{mainstab2}, \ref{mainstab3} and \ref{instap}. 
Our orbital stability results will follow by some strict 
relationship between different minimization problems.

\bdf\label{l2}
Given $\gamma>0$, let us consider the minimization problems
\begin{align}
\label{elledueI}
\min_{\M}\I(U) &=m_{\ga},\quad
\\
\label{elledueE}\min_{\M}\E(U)&=c_{\ga},
\end{align}
where $\M=\left\{ U\in\H : \|U\|_{2,\om}^{2}=\ga\right\}$. 
Moreover, we denote with $\A$ the set of the solution of 
problem \eqref{elledueE}.
\edf
\bos\label{elledueIE}
Notice that, solving problem \eqref{elledueI} is equivalent 
to solve problem \eqref{elledueE}, since for every 
$V\in\M$ we have $\I(V)=\E(V)+\gamma/2$.
\eos
\bos\label{reali}
It results
\begin{equation*}
\A=\{(e^{i\theta_1}u_1,e^{i\theta_2}u_2),\,
\theta_{j}\in \R,\,(u_{1},u_{2})\in H^{1}(\R^{n},\R^{2})\,
\text{solves}\,\eqref{elledueEr}\},
\end{equation*}
where
\beq\label{elledueEr}
\sigma_{\R} =\inf \left\{\E(V): V\in H^1(\R^n;\R^{2}),\;
\om_{1}\|v_{1}\|^{2}_{2}+\om_{2}\|v_{2}\|^{2}_{2}=\gamma \right\}.
\eeq
Indeed, if we consider the minimization problems \eqref{elledueE}
we have that $c_{\gamma}=\sigma_{\R}$ and if $U=(u_{1},u_{2})$ 
solves \eqref{elledueE} there exists $\theta_{j}\in \R$
such that  $u_{j}=e^{i\theta_{j}}|u_{j}|$ for $j=1,2$.
For more details, see Remark 3.12 of \cite{mopesq}.
\eos

\bos
The conservation laws of the problem suggest that orbital stability has
to be studied by using problem \eqref{elledueI} which can be solved only 
for $p<1+2/n$, as for $p>1+2/n$, $\E$ (and then $\I$) is not bounded from below on $\M$. We will prove our stability result using problem \eqref{cne}
which has a solution for every $p<n/(n-2)$. At the same time we cannot expect to have a stability result for every $p<n/(n-2)$, as Theorem \ref{instap}
shows. This aspect is clarified in the following results where we show that problems \eqref{cne} and \eqref{elledueE} (and then $\eqref{elledueI}$)
are equivalent for $p<1+2/n$. Indeed, in this range of exponents 
we can construct a bijective correspondence between the negative 
critical values of $\E$ on $\M$ and the critical values of $\I$ on 
$\nh$. While, for $p>1+2/n$ we cannot derive this map between 
these critical values.
\\
This suggests that for $p<1+2/n$ the Nehari manifold and $\M$ 
have the same tangent planes, while when $p>1+2/n$ the tangent 
planes are different, so that a minimum point on $\nh$ would 
probably give rise to a different critical point on $\M$. 
This point is crucial in proving orbital stability properties, 
since the conservation laws show that the dynamical analysis 
has to be performed on $\M$.
\eos

In proving many of the results of this section we will make 
use of the following lemma the proof of which is straightforward.
\begin{lemma}\label{ulamu}
For any $u\in H^{1}(\R^{n},\C)$ and for any positive real numbers 
$\lambda, \mu$, we can define the scaling 
$u^{\mu,\lambda}(x)=\mu u(\lambda x)$
such that the following equalities hold.
\beq\label{scaling}
\|u^{\mu,\la}\|_{2}^{2}=\mu^{2}\la^{-n}\|u\|_{2}^{2},\quad
\|\nabla u^{\mu,\la}\|_{2}^{2}=\mu^{2}\la^{2-n}\|\nabla u\|_{2}^{2},
\quad \|u^{\mu,\la}\|_{2p}^{2p}=\mu^{2p}\la^{-n}\|u\|_{2p}^{2p}.
\eeq
\end{lemma}

\bdim
The proof is an immediate consequence of a change of variables.
\edim

First, we want to show the equivalence between problems 
\eqref{cne} and \eqref{elledueE} and \eqref{elledueI}. 
In order to do this, let us define the sets
\bdm
K_{\mathcal E}\!=\!\left\{c<0\,:\,\exists \,U\in \M\,
:\, \E(U)=c,\,\nabla_{\!\!\M}\E(U)=0 \right\},
\qquad
{\tilde K}_{\mathcal E}\!=\!
\left\{U\in \M\,:\,\nabla_{\!\!\M}\E(U)=0,\; \E(U)<0\right\},
\edm
\beq\label{defkI} 
K_{\mathcal I}\!=\!\left\{m\in\R\,:\,\exists \,U\in \mathcal{N}\,
:\,\, \I(U)=m,\,\,\I'(U)=0 \right\},
\qquad
{\tilde K}_{\mathcal I}\!=\!\left\{U\in \mathcal{N}\,:\,
\I'(U)=0 \right\},
\eeq
where we have denoted with $\nabla_{\!\!\M}\E $ the tangential 
derivative of $\E$ on $\M$. The following result holds.

\bte\label{critici}
Assume \eqref{pzero}. For any $\bt,\om_1,\om_2>0$ the following 
conclusions hold: \\
$a)$ there exists a bijective correspondence between the sets
${\tilde K}_{\mathcal E}$ and ${\tilde K}_{\mathcal I}$,\\
$b)$ there exists a  bijective map 
$T:K_{\mathcal I}\to K_{\mathcal E}$ given by
\beq\label{mec}
\begin{array}{l}
T(m)=-\left[\dfrac{1}{p-1}-\dfrac{n}{2}\right] 
\left[\dfrac{\ga}{2p'-n}\right]^{\frac{2p'-n}{2/(p-1)-n}}
\left[\dfrac{1}{m }\right]^{\frac{2}{2/(p-1)-n}},
\end{array}
\eeq
where $p'=p/(p-1)$ stands for the conjugate exponent of $p$.
\ete

\bdim In order to prove assertion {\it a)}, take 
$V=(v_{1},v_{2})\in\M$ such that $V$ satisfies
\beq\label{evincolato}
\langle \E'(V), V \rangle= -\nu\ga,\quad
\E(V)=c<0.
\eeq
Then, using that $\langle F'(V),V\rangle=2p F (V)$,
it follows
\bdm
-\nu\ga-2c=\langle \E'(V),V \rangle-2\E(V)=2F(V)(1-p)<0,
\edm
showing that $\nu$ is a positive real number. Therefore
it is well defined the map $T^{\mu,\la}:\M\to {\mathcal N}$ 
\bdm
T^{\mu,\la}(V)=V^{\mu,\la},
\edm
where $\mu,\,\la$, are given by
\beq\label{parameters}
\mu=\nu^{-1/2(p-1)},\quad  \la=\nu^{-1/2}.
\eeq
Using this and \eqref{evincolato} one obtains that $V^{\mu,\la}$ 
solves \eqref{ellittico}, so that $T^{\mu,\la}(V)$ belongs to 
${\tilde K}_{\mathcal I}$. Vice-versa if 
$U\in {\tilde K}_{\mathcal I}$ let us take $\nu>0$ such that
\bdm
\nu^{1/(p-1)-n/2}=\frac{\gamma}{\|U\|_{2,\om}^{2}}
\edm
and $\la,\,\mu>0$ given by
\bdm
\la=\nu^{1/2},\qquad \mu=\nu^{1/2(p-1)},
\edm
so that $U^{\mu,\la}$ belongs to $\M$, $\nabla_{\M}\E(U)=0$.
This shows that $ (T^{\mu,\la})^{-1}=T^{1/\mu,1/\la}$.
\\
In order to prove assertion {\it b)}, note first that any 
$m\in K_{\mathcal I}$ is positive. Indeed, since there exists 
$U\in \mathcal{N}$ such that $\I(U)=m$ and $\I'(U)=0$, it follows
\bdm
m=\I(U)-\dfrac{1}{2p}\langle\I'(U),U\rangle
=\dfrac{1}{2}\left(1-\dfrac{1}{p}\right) \|U\|^{2}_{\H}>0,
\edm
so that $T$ is a well defined and injective map.
Let us first show that if $c\in K_{\mathcal E}\cap\R^{-}$,
then $c=T(m)$.

Indeed take $V\in\M$ corresponding to such $c$, and take
$T^{\mu,\la}(V)=V^{\mu,\la}$. Recalling Pohozaev identity 
(see $(5.9)$ in \cite{mmp}) and since 
$V^{\mu,\la}\in\mathcal{N}$ we get
\bdm
\left(\dfrac{1}{2}-\dfrac{1}{2p}\right)
\left(\|\nabla V^{\mu,\la}\|_{2}^{2}
+\|V^{\mu,\la} \|_{2,\om}^{2}\right)=m, 
\edm 
\bdm
(n-2)\|\nabla V^{\mu,\la}\|_2^2 +n\|V^{\mu,\la}\|_{2,\om}^{2}
=\dfrac{n}{p} \left(\|\nabla V^{\mu,\la}\|_{2}^{2}
+\|V^{\mu,\la} \|_{2,\om}^{2}\right),
\edm
where $m=\I(V^{\mu,\la})$. We derive
\beq\label{poho}
\|\nabla V^{\mu,\la}\|_{2}^{2}
= n m,\qquad F (V^{\mu,\la})
=\dfrac{m}{p-1},\qquad \|V^{\mu,\la} \|_{2,\om}^{2}
=\left(\dfrac{2p}{p-1}-n\right) m.
\eeq
Using \eqref{scaling} we have that
\beq\label{scaling2}
\begin{array}{l}
\dys\frac{\mu^{2}}{\la^{n-2}}\|\nabla V\|_{2}^{2}= n m, \quad
\dys \frac{\mu^{2p}}{\la^{n}}F (V) =\dfrac{m}{p-1}, \quad
\dys \frac{\mu^{2}}{\la^{n}}\|V\|_{2,\om}^{2}
=\left(\dfrac{2p}{p-1}-n\right) m.
\end{array}
\eeq

Since $V$ is in $\M$, \eqref{parameters} yields
\bdm
\ga=\|V\|_{2,\om}^{2} =\nu^{1/(p-1)-n/2}
\left(\dfrac{2p}{p-1}-n\right) m,
\edm
this and \eqref{scaling2} give
\bdm
\begin{array}{l}
\|\nabla V\|_{2}^{2}=n\left(\dfrac{\ga}{2p'-n}
\right)^{1+\frac{2(p-1)}{2-n(p-1)}}
\left(\dfrac{1}{m}\right)^{\frac{2(p-1)}{2-n(p-1)}},
\\ \\
F(V)= \dfrac{1}{p-1}\left(\dfrac{\ga}{2p'-n}
\right)^{1+\frac{2(p-1)}{2-n(p-1)}}
\left(\dfrac{1}{m}\right)^{\frac{2(p-1)}{2-n(p-1)}}.
\end{array}
\edm

All the above calculations imply that, if $c$ is a negative 
constrained critical value of $\E$ on $\M$ and $m$ is the 
corresponding critical value of $\I$, than $c$ is given by 
\eqref{mec}.

In order to show that $T^{-1}$ is surjective let us take
$m$ in $K_{\mathcal I}$ and the corresponding $U$ that satisfies
the conditions in \eqref{defkI}. For any $\nu>0$ we can define
$$
\lambda=\nu^{1/2},\qquad \qquad \mu=\nu^{1/2(p-1)}
$$
and consider $U^{\mu,\lambda}$. Using \eqref{poho} and 
\eqref{scaling} and requiring that $U^{\mu,\lambda}\in\M$ 
imply that $\nu$ is related to $\gamma$ by the expression
$$
\nu^{1/(p-1)-n/2}=\frac{\gamma}m\left[\frac{1}{2p'-n}\right].
$$
Moreover, since $U$ is a free critical point of $\I$ we obtain 
that $U^{\mu,\lambda}$ is a constrained critical point of $\E$
with Lagrange multipliers equal to $\nu$. In order to conclude 
the proof we have to impose that $\E(U^{\mu,\lambda})=c$.
>From conditions \eqref{poho}, \eqref{scaling} and from the 
definition of $K_{\mathcal I}$ it follows that $c,\,m$ 
and $\nu$ satisfy
$$
c=m\left[\frac{n}2-\frac1{p-1}\right]\nu^{p'-n/2}
$$
and substituting the value of $\nu$ in dependence of $\gamma$ 
implies that $m=T^{-1}(c)$.
\edim

\begin{corollary}\label{equiground}
There exists a bijective correspondence between the sets 
$\G$ and $\A$.
\end{corollary}

\bdim
Let $V\in\A$ and take $T^{\mu,\la}(V)$; Theorem \ref{critici}
implies that $T^{\mu,\la}(V)$ is a critical point of $\I$, so that
we only have to show that 
$$
\I(T^{\mu,\la}(V))=m=\mnehari.
$$ 
Indeed, suppose by contradiction that $m>\mnehari$. 
In \cite{mmp} it is proved that 
$\mnehari$ is achieved by a vector $U$, then $U^{1/\mu,1/\la}$,
with $\mu,\la$ as in \eqref{parameters}, belongs to $\M$ and  
gives a negative critical value $c$ given by \eqref{mec}.
Since $\mnehari<m$ we get $c<\cdue$ which is a contradiction, 
so that the claim is true.  
\edim
Using the preceding result and Theorem $2.1$ in \cite{mmp}
we can prove the following statement.

\bte\label{groundstate}
Assume \eqref{pzero}. For any $\bt,\om_1,\om_2>0$, there
exists a solution of the minimization problems \eqref{elledueI}, \eqref{elledueE}.
\ete

\bdim 
As observed in Remark \ref{elledueIE} problems \eqref{elledueI}, \eqref{elledueE} are equivalent, so it is enough to show that 
\eqref{elledueE} is solved.
\\
By using a Gagliardo-Nirenberg type inequality for systems 
(see \cite{fm} equation $(9)$), we get that the following 
inequality holds for any $U\in\M$ 
\bdm
\E(U)\geq\dfrac{1}{2}\|\nabla U\|_{2}^{n(p-1)}
\left[\|\nabla U\|_{2}^{2-n(p-1)} 
-\dfrac{C_{\om,\bt}}{p}\ga^{p-(p-1)n/2}\right],
\edm
so that $\E$ is bounded from below if and only if \eqref{pzero} 
holds. Moreover, note that the infimum $\cdue$ in \eqref{elledueE} 
is negative. 
Indeed, we impose $\la^n=\mu^2$ so that, for any $U\in\M$,
$U^{\mu,\la}= \left( u_{1}^{\mu,\la},u_{2}^{\mu,\la}\right)$ 
still belongs to $\M$. By \eqref{scaling} we derive the real 
function $h$ defined by
\bdm
h(\la)=\E (U^{\mu,\la})= \dfrac{\la^2}{2} \|\nabla U\|_{2}^{2}
-\la^{n(p-1)}F (U),
\edm
and from condition \eqref{pzero} it follows that there exists a 
$\la_0=\la_0(U)>0$  such that for any $\la\in(0,\la_0)$ $h(\la)$ 
is negative and  this shows the claim. Then, Theorem \ref{critici}
implies that there exists $m\in K_{\mathcal I}$ such that 
$\cdue=T(m)$. Finally, since in \cite{mmp} it is proved that 
$\mnehari$ is achieved, using Corollary \ref{equiground} we 
obtain the conclusion.
\edim

\bos\label{gamma}
It is easy to see that every $U$ in $\G$ satisfies
\bdm
\|U\|_{2,\om}^{2}=m_{\nh}\left(\frac{2p}{p-1}-n\right),
\edm
thanks to the regularity properties of $U$ and to Pohozaev identity.
\eos

\bte\label{muguali}
Assume \eqref{pzero} and let $\gamma_{0}$ be fixed as
\beq\label{defgamma}
\gamma_{0}=m_{\nh}\left(\frac{2p}{p-1}-n\right).
\eeq
Then $m_{\nh}=m_{\gamma_{0}}$.
\ete

\bdim From the definition of $\ga_{0}$ immediately follows that 
$m_{\nh}\geq m_{\ga_{0}}$. In order to show that the 
equality is achieved, we only have to observe that
$m_{\ga_{0}}=c_{\ga_{0}}+\ga_{0}/2=T(m_{\nh})+\ga_{0}/2$.
Using the definition of $T$ joint with \eqref{defgamma} yields the conclusion.
\edim

\bos
Consider the minimization problems \eqref{cne} and
\beq\label{cner}
m_{\nh_{\R}} :=\inf\{\I(W) \,:\,W\in H^{1}(\R^{n},\R^{2})\setminus\{0\}
\,:\,\langle \I'(W),W\rangle =0\,\},
\eeq
it results that $m_{\nh}=m_{\nh_{\R}}$. Indeed, $m_{\nh}\leq m_{\nh_{\R}}$. 
Moreover, if $U\in\H,\,U=(u_{1},u_{2})$ is a solution of \eqref{cne} then
$\|U\|_{2,\om}=\gamma_{0}$, where $\gamma_{0}$ is defined in 
\eqref{defgamma}, so that $U\in {\mathcal M}_{\gamma_{0}}$ and Theorem
\ref{muguali} implies that $\E(U)=c_{\gamma_{0}}$. Then, from
Remark \ref{reali} we deduce that there exist $\theta_{1},\,\theta_{2}$
such that
$$
U=(u_{1},u_{2})=(e^{i\theta_{1}}|u_{1}|,e^{i\theta_{2}}|u_{2}|), 
\qquad \E(|u_{1}|,|u_{2}|)=c_{\gamma_{0}}.
$$
Finally, Theorem \ref{muguali} gives that $\I(|u_{1}|,|u_{2}|)=m_{\nh}$
and, since $|\nabla u_{i}|=|\nabla |u_{i}||$ it results that $(|u_{1}|,|u_{2}|)
\in \nh$ is a solution of \eqref{cner}, that is $m_{\nh}$ is achieved on 
a vector with real valued components. Furthermore, we have shown that
$$
\G=\{(e^{i\theta_{1}}u_{1},e^{i\theta_{2}}u_{2}),\,\theta_{j}\in\R,
\,(u_{1},u_{2})\in H^{1}(\R^{n},\R^{2})\,\text{solves}\,\eqref{cner}\}.
$$
\eos

In order to prove the instability result Theorem \ref{instap} another
variational characterization of a ground state solution will be useful.
Let us define the functional 
\beq\label{erre}
\mathcal{R}(U)=\|\nabla U\|_{2}^{2}-n(p-1) F(U)
\eeq
and the infimum
$$
m_{\P}=\inf_{\P}\I\;\text{ where }\;\P=\{U\in \H\,:\,\mathcal{R}(U)=0\}.
$$
The following results hold.

\bpr\label{infs}
Assume that $p>1+2/n$, then the following conclusions hold:\\
$a)$ $\P$ is a natural constraint for $\I$;\\
$b)$ $m_{\P}=\mnehari$.
\epr

\bdim In order to prove $a)$ let us consider $U$ a constrained 
critical point of $\I$ on $\P$, then there exists $\la\in\R$ 
such that the following identities are satisfied
\beq\label{uno}
(1-2\la)\|\nabla U\|_{2}^{2}+\|U\|_{2,\om}^{2} =2p[\la n(1-p)+1]F(U),
\eeq
\beq\label{due}
(1-2\la)\left(\dfrac{n}{2}-1\right)\|\nabla U\|_{2}^{2}+\dfrac{n}{2}
\|U\|_{2,\om}^{2} =n[\la n(1-p)+1] F(U),
\eeq
\beq\label{tre}
\|\nabla U\|_{2}^{2} =n(p-1) F(U).
\eeq
Hence, using \eqref{tre} in \eqref{uno} we get
\bdm
\|U\|_{2,\om}^{2}=\left[2\la(1-p)+\dfrac{2p}{n(p-1)}-1\right]
\|\nabla U\|_{2}^{2},
\edm
and using this and \eqref{tre} in \eqref{due}, and taking into 
account that $p>1+2/n$, we obtain that $\la=0$, so that $U$ is 
a free critical point of $\I$. 

In order to prove $b)$ take a minimum point $U$ of $\I$ in $\P$; 
from $a)$ it follows that then $U$ belongs to $\mathcal{N}$ so 
that $m_{\P}\geq \mnehari;$ viceversa if $V$ is a minimum point 
of $\I$ in $\mathcal{N}$ then $V$ is a free critical point of 
$\I$ and Pohozaev identity implies that $V\in\P$ so that 
$m_{\P}\leq \mnehari$, yielding the conclusion.
\edim

\ble\label{glambda}
Assume that $p>1+2/n$.
For any $U\neq(0,0)$ let us consider 
$U^{\la^{n/2},\la}=\left(u_{1}^{\la^{n/2},\la},
u_{2}^{\la^{n/2},\la}\right)$, for $u^{\mu,\la}$ defined 
in Lemma \ref{ulamu}. The following conclusions hold: 
\\
$a)$ there exists a unique $\la_{*}=\la_{*}(U)$ such that 
$U^{\la_{*}^{n/2},\la_{*}}$ belong to $\P$,
\\
$b)$ the function $g(\la)=\I\left(U^{\la^{n/2},\la}\right)$
has its unique maximum point in $\la=\la_{*}$,
\\
$c)$ $\la_{*}<1$ if and only if ${\mathcal R}(U)<0$ and 
$\la_{*}=1$ if and only if $\mathcal{R}(U)=0$,
\\
$d)$ the function $g(\la)$ is concave on 
$\left(\la_{*},+\infty\right)$.
\ele

\bdim $a)$: for any $\la>0$ it holds
\bdm
\mathcal{R}\left(U^{\la^{n/2},\la}\right)=\la^{2}\|\nabla U\|^{2}_{2}
-\la^{n(p-1)}n(p-1)F(U),
\edm
then there is a unique 
\bdm
\la_{*}(U)=\la_{*}=\left[\dfrac{\|\nabla U\|^{2}_{2}}{n(p-1)F(U)}
\right]^{1/[n(p-1)-2]}
\edm
such that $\mathcal{R}\left(U^{\la_{*}^{n/2},\la_{*}}\right)=0$.
Computing the first derivative of $g(\la)$ $b)$ is proved.
Since $p>1+2/n$, $c)$ easily follows.
$d)$ immediately follows from writing the second derivative
of the function $g$.
\edim
\ble\label{RI}
For any $U\in\H$ with $\mathcal{R}(U)<0$ it results
$$
\mathcal{R}(U)\leq \I(U)-\mnehari.
$$
\ele

\bdim Conclusion {\it d)} in Lemma \ref{glambda} implies that
$$
g(1)\geq g(\la_{*})+g'(1)(1-\la_{*}).
$$
Direct computation yields
\bdm
\I(U)=g(1)\geq g(\la_{*})+g'(1)(1-\la_{*})
=g(\la_{*})+ \mathcal{R}(U) (1-\la_{*})\geq 
\I\left(U^{\la_{*}^{n/2},\la_{*}}\right)+ \mathcal{R}(U),
\edm
where in the last inequality we have used that $\mathcal{R}(U)<0$.
Recalling that $U^{\la_{*}^{n/2},\la_{*}}\in\P$ and
applying conclusion $b)$ of Proposition \ref{infs}
complete the proof.
\edim


\subsection{Bound states}

It is well known (see Section 8 in \cite{caz}) that $\zbo$ defined in 
\eqref{defzbo} can be characterized as the solution of 
the following constrained minimization problem
\beq\label{miniequa}
\E_{1}(u)=c_{\delta}=\min_{\Md}\E_{1} \quad\text{where}
\quad \Md=\{u\in H^{1}(\R^{n},\R)\,:\,\|u\|_{2}^{2}=\delta\};
\eeq
where the functional $\E_{1}:H^{1}(\R^{n})\to \R$ is defined by
\beq\label{defE1}
\E_{1}(u)=\frac12\|\nabla u\|_{2}^{2}-\frac{\beta+1}{2p}\|u\|_{2p}^{2p}.
\eeq
and when we prescribe 
$$
\delta=\delta(\om)=\frac{\om^{\frac1{p-1}-\frac{n}2}}{
(\beta+1)^{\frac1{p-1}}}\|z^{1}_{0}\|_{2}^{2}.\quad 
$$ 
Otherwise, $\zbo$ can be equivalently obtained as
the solution of the minimization problem
\beq\label{nehari1}
\I_{1}(u)=m_{1}=\min_{\mathcal{N}_{1}}\I_{1}\quad
\hbox{ where }\quad 
\mathcal{N}_{1}=\left\{u\in H^{1}(\R^{n},\R),
\,u\not\equiv 0:\, \langle \I_{1}'(u),u\rangle=0\right\};
\eeq
where the functional $\I_{1}:H^{1}(\R^{n},\R)\to \R$ is defined by
$$
\I_{1}(u)=\dfrac{1}{2}\|\nabla u\|^{2}_{2}+\frac{\om}2\|u\|_{2}^{2}
-\frac{\beta+1}{2p}\|u\|_{2p}^{2p}.
$$
The following result is the starting point in proving Theorem \ref{carabound}.
\bpr\label{moduliN}
Assume that $\om_{1}=\om_{2}=\om$. If\, 
$U$ solves \eqref{cnedue} it results $|u_{1}|=|u_{2}|$ 
almost everywhere.
\epr
\bdim
Consider the variational characterization $\zbo$ as  the solution of \eqref{nehari1}, the 
vector $Z=(\zbo,\zbo)$ belongs to $\nh_{2}$, so that 
$$
2m_{1}= 2\I_{1}(z)=\I(Z)\geq m_{2}.
$$
Let now $U$ be a solution of \eqref{cnedue}, then, 
Young inequality yields
$$
m_{2}=\I(U)\geq \I_{1}(u_{1})+\I_{2}(u_{2})\geq 2m_{1}\geq m_{2}
$$
showing that 
\beq\label{equaI}
\I(U)=\I_{1}(u_{1})+\I_{1}(u_{2}),\quad \text{and}\quad m_{2}=2m_{1}.
\eeq 
Writing down this equality we get
$$
\frac{\beta}{p}\int |u_{1}u_{2}|^{p}=\frac{\beta}{2p}
\int \left(|u_{1}|^{2p}+|u_{2}|^{2p}\right),
$$
that is
$$
\int\left(|u_{1}|^{p}-|u_{2}|^{p}\right)^{2}=0,
$$
giving the conclusion.
\edim
{\bf Proof of Theorem \ref{carabound}.}\quad 
Let $U=(u_{1},u_{2})\in \B$, from  Proposition \ref{moduliN} we derive that
$$
\|\nabla u_{1}\|_{2}^{2}+\om\|u_{1}\|_{2}^{2}=\|u_{1}\|_{2p}^{2p}+
\bt\|u_{1}u_{2}\|_{p}^{p}=(\bt+1)\|u_{1}\|_{2p}^{2p},
$$
that is $u_{1}\in \mathcal{N}_{1}$, so that 
\beq\label{disI}
\I_{1}(u_{1})\geq \I_{1}(\zbo).
\eeq 
On the other hand, from Proposition \ref{moduliN} and \eqref{equaI} 
it follows that $\I(U)=2\I_{1}(u_{1})$
and recalling that $Z=(\zbo,\zbo)\in \nh_{2}$ we derive
$$
2\I_{1}(u_{1})=\I(U)\leq \I(Z)=2\I_{1}(\zbo).
$$
This, \eqref{disI} and \eqref{nehari1} imply that $u_{1}=e^{i\theta_{1}}\zbo(\cdot-y_{1})$ for some $\theta_{1}\in \R$ and $y_{1}\in \R^{n}$. 
The same argument for $u_{2}$ gives $u_{2}=e^{i\theta_{2}}
\zbo(\cdot-y_{2})$, and Proposition
\ref{moduliN} yields $y_{1}=y_{2}$. 
\edim
As we did for ground states we want to investigate the connection of
Problem \eqref{cnedue} with a minimization of $\E$ under  
suitable constraints. Since we are now considering vectors with both nontrivial components we are naturally lead to study the following problem for $\E$.
\beq\label{cmdue}
\E(U):=c_{(\delta_{1},\delta_{2})}=\inf_{\Mdue} \E \quad\text{where}\quad
\Mdue = \{(u_{1},u_{2})\in\H \; : \;
\|u_{1}\|^{2}_{2}=\delta_{1}, \; \|u_{2}\|^{2}_{2}=\delta_{2}\}
\eeq
where $\delta_{1},\delta_{2}$ are positive real numbers. 


If $\delta_{1}\neq \delta_{2}$ we do not know how to solve 
problem \eqref{cmdue}, and, to our knowledge, it 
is an open problem to prove orbital stability property of solution of \eqref{cmdue} in this general case.
\\
Therefore, we will focus our attention to the case $\delta_{1}=\delta_{2}=
\dt$.
In this case we have the following minimization problem
\beq\label{ohta}
\E(U)=c_{(\delta,\delta)}=\inf_{\Mdd}\E \quad\text{where}
\quad \Mdd=\left\{ (u_{1},u_{2})\in\H \; : \; 
\|u_{1}\|_{2}^{2}=\|u_{2}\|_{2}^{2}=\delta\right\}.
\eeq
As we did for the ground state solutions, investigating the relation 
between Problems \eqref{ohta} and \eqref{cnedue} in the case 
$\om_{1}=\om_{2}=\om$ naturally lead us to choose $\delta=\delta(\om)$ 
such that every solution of \eqref{ohta} give rise a solution of 
Problem \eqref{ellittico}. With this choice we end up with the same 
characterization of the sets ${\mathcal B}$ (given in Theorem 
\ref{carabound}) and ${\mathcal A}_{(\delta(\om),\delta(\om))}$ 
set of solutions of \eqref{ohta}, as the next result shows.
\\
The following result can be proved using the same arguments of \cite{oh}
for the case $p=2$ and $\bt=1$. We include some details for clearness.
\bpr\label{caraA}
It results
$$ 
\A_{(\delta(\om),\delta(\om))}=\left\{
(e^{i\theta_{1}} \zbo(\cdot-y),e^{i\theta_{2}} 
\zbo(\cdot-y)),\;\theta_{i}\in\R,\, y\in\R^{n}\right\}.
$$
\epr
\bdim
Taking $U=(u_{1},u_{2})\in \A_{(\delta(\om),\delta(\om))}$, 
using  the variational characterization $\zbo$ as  the solution of \eqref{miniequa}, and  arguing as in the proof of Proposition \ref{moduliN}
yield
\beq\label{moduli}
|u_{1}|=|u_{2}|, \qquad \E(U)=\E_{1}(u_{1})+\E_{1}(u_{2}),\qquad
c_{(\delta,\delta)}=2c_{\delta}
\eeq
Moreover, $\|u_{i}\|_{2}^{2}=\delta(\om)$, so that 
$u_{i}\in \Md$ and it results $\E_{1}(u_{i})\geq c_{\delta}$.
Then \eqref{moduli} yields $\E_{1}(u_{i})=c_{\delta}$.
Thus, $u_{i}$ are solutions of the minimization problem 
\eqref{miniequa}, and we can use the uniqueness result in \cite{k}
and Theorem II.1 in \cite{cl} to deduce that there exist
$\theta_{1},\,\theta_{2},\,y_{1},\,y_{2}$ such that 
$u_{1}=e^{i\theta_{1}}\zbo(\cdot-y_{1})$,  
$u_{2}=e^{i\theta_{2}}\zbo(\cdot-y_{2})$. 
Moreover, \eqref{moduli} implies $y_{1}=y_{2}$. 
\edim

\section{Proofs of the Main Results}\label{proofs}

In this section we will prove the main results concerning 
the stability (or instability) of the standing waves. In 
particular in the following subsections we show, in the subcritical 
case $1<p<1+2/n$,  the orbital stability of the sets $\G$ and $\Se$,
and also of the set $\B$ for $\om_{1}=\om_{2}$.
Finally in subsection \ref{subcrit} we prove that for $p>1+2/n$ 
the sets $\G,\, \Se$ and $\B$ (for $\om_{1}=\om_{2})$ 
are unstable and for $p=1+2/n$ the instability holds for every 
bound state. The proofs of Theorems \ref{mainstab}, \ref{mainstab2} 
\ref{instap} and \ref{instapB} follow the arguments  of \cite{cl,caz} 
for the single equation, while the proof of Theorem 
\ref{mainstab3} follows the arguments of \cite{oh}.

\subsection{Stability of the ground state standing waves}

{\bf Proof of Theorem \ref{mainstab}.}\quad 
Let us argue by contradiction, and suppose that there exist 
$\eps_{0}>0$, $\{t_{k}\}\subset\R$ and a sequence of initial 
data $\{\Phi^{0}_{k}\}\subset\H$ such that
\beq\label{ipodati}
\lim_{k\to \infty}\inf_{U\in\G}\|\Phi^{0}_{k}-U\|_{\H}=0
\eeq
and the corresponding sequence of solution $\{\Phi_{k}\}$ of
Problem \eqref{schr} satisfies
\beq\label{iposol}
\inf_{U=(u_{1},u_{2})\in \G}\|\Phi_{k}(\cdot, t_{k})
-(u_{1},u_{2})\|
\geq \eps_{0}.
\eeq
Condition \eqref{ipodati}, definitions \ref{ground}, 
\eqref{defgamma},  Theorem \ref{muguali} and the 
continuity properties of the functional $\I$ yield
$$
\I(\Phi^{0}_{k})\to m_{\gamma_{0}}\qquad   
\|\Phi^{0}_{k}\|^{2}_{2,\om}=\gamma_{0}+\text{\small o}(1).
$$
Let us denote $\Psi_{k}(x)=\Phi_{k}(x,t_{k})$,
then, conservation laws \eqref{mass}, \eqref{energy} imply
that  
$$
\I(\Psi_{k})\to m_{\gamma_{0}},\qquad 
\|\Psi_{k}\|^{2}_{2,\om}=\gamma_{0}+\text{\small o}(1).
$$
Following the arguments of \cite{mmp}, we use the Ekeland variational 
principle to obtain a new minimizing sequence $\tilde{\Psi}_{k}$ 
which is also a Palais-Smale sequence for $\I$ and 
we find $\tilde{\Psi}$ such that
\bdm
\I(\tilde{\Psi})=m_{\gamma_{0}},\quad \I'(\tilde{\Psi})=0.
\edm
Then 
\bdm
\dfrac{1}{2}\left(1-\dfrac{1}{p}\right)\|\tilde{\Psi}\|^{2}_{\H}
=m_{\gamma_{0}}=\I(\tilde{\Psi}_{k})-\dfrac{1}{2p}
\langle \I'(\tilde{\Psi}_{k}),\tilde{\Psi}_{k} \rangle+o(1)
=\dfrac{1}{2}\left(1-\dfrac{1}{p}\right)\|\tilde{\Psi}_{k}\|^{2}_{\H}
+o(1)
\edm
showing that
\bdm
\|\tilde{\Psi}_{k}\|^{2}_{\H}\to\|\tilde{\Psi}\|^{2}_{\H}
\edm
and, since $\H$ is an Hilbert space, we have the strong
convergence in $\H$. By the choice of $\tilde{\Psi}_{k}$ 
we derive that also $\Psi_{k}\to\tilde{\Psi}$ strongly in $\H$,
which is an evident contradiction with \eqref{iposol}.
\edim
\bos
In the proof of the previous result it was crucial to show that 
the sequence $\tilde{\Psi}_{k}$ strongly converges in $\H$, 
to get the desired contradiction. In other words, in proving orbital
stability results we made use of conservation laws, minimization property and compactness of the minimizing sequence.
\eos
\subsection{Stability of bound state standing waves}

{\bf Proof of Theorem \ref{mainstab2}.}\quad
Let us argue for the set of scalar solution with the second component 
equal to zero, the other case can be handled analogously.
The conclusion can be obtained arguing as in the previous 
Theorem assuming that there exist $\eps_{0}>0$, 
$\{t_{k}\}\subset\R$ and a sequence of initial data 
$\{\Phi^{0}_{k}\}\subset\H$ such that
\bdm
\lim_{k\to \infty}\inf_{\theta\in\R,y\in \rn}
\|\Phi^{0}_{k}-(e^{i\theta}z^{\om}_{0}(\cdot-y),0)\|_{\H}=0
\edm
and the corresponding sequence of solution $\{\Phi_{k}\}$ 
of Problem \eqref{schr} satisfies
\beq\label{absurd}
\inf_{\theta\in\R,y\in \rn}\|\Phi_{k}(\cdot, t_{k})
-(e^{i\theta}z^{\om}_{0}(\cdot-y),0)\|
\geq \eps_{0}
\eeq
As before, via conservation laws, 
$\Psi_{k}(x)=\Phi_{k}(x,t_{k})$ satisfies
\beq\label{consescal}
\E(\Psi_{k}(t))\to \E(z^{\om}_{0},0)=c_{(\delta_{0},0)},\qquad 
\|\psi_{k,1}\|_{2}=\|z^{\om}_{0}\|_{2}+\text{\small o}(1),\;
\|\psi_{k,2}\|_{2}=\text{\small o}(1),
\eeq
where
$$
\E(z^{\om}_{0},0)=c_{(\delta_{0},0)}=\min_{\mathcal{M}_{(\delta_{0},0)}}\E\,,\quad
\quad \mathcal{M}_{(\delta_{0},0)}=\{U\in\H\;:\;
\|u_{1}\|_{2}=\|z^{\om}_{0}\|^{2}_{2}=\delta_{0},\;\|u_{2}\|^{2}_{2}=0\}.
$$
>From Gagliardo-Nirenberg inequality we deduce that $\Psi_{k}$ 
is bounded in $\H$, then interpolation inequality, joint with \eqref{consescal}, implies that 
\beq\label{convforte}
\psi_{k,2}\to 0 \quad \text{ strongly in $L^{2p}$}. 
\eeq
Therefore, Holder inequality yields
$$
\int |\psi_{k,1}|^{p}|\psi_{k,2}|^{p}\leq 
\|\psi_{k,1}\|_{2p}^{p}\|\psi_{k,2}\|_{2p}^{p}\to 0.
$$
This and \eqref{convforte} yield
\beq\label{gradient}
\E(\Psi_{k})=\E_{0}(\psi_{k,1})
+\frac12\|\nabla \psi_{k,2} \|_{2}^{2}+\text{\small o}(1),
\eeq
where 
$$
\E_{0}(u)=\frac12\|\nabla u\|^{2}_{2}-\frac1{2p}\|u\|_{2p}^{2p}.
$$
So that we are lead to 
$$
\E_{0}(\psi_{k,1})\leq \E(\Psi_{k})+\text{\small o}(1)=c_{(\delta_{0},0)}
+\text{\small o}(1).
$$
Consider the scalar minimization problem
$$
\E_{0}(z^{\om}_{0})=c_{\delta_{0}}=\min_{\mathcal{M}_{\delta_{0}}}\E_{0}\qquad 
\mathcal{M}_{\delta_{0}}=\{u\in H^{1}(\R^{n},\R)\,:\,\|u\|_{2}^{2}=
\|z^{\om}_{0}\|_{2}^{2}\}.
$$
It is easy to verify that 
\beq\label{equac}
c_{\delta_{0}}=c_{(\delta_{0},0)},
\eeq
so that
\beq\label{scalarc}
\E_{0}(\psi_{k,1})\leq\E_{0}(z^{\om}_{0})+\text{\small o}(1)
=c_{\delta_{0}}+\text{\small o}(1),
\qquad 
\|\psi_{k,1}\|_{2}^{2}\to \|z^{\om}_{0}\|_{2}^{2}.
\eeq
By using Gagliardo-Nirenberg type inequality and
arguing as in the proof of Theorem \ref{groundstate}, we obtain that $c_{(\delta_{0},0)}<0$, so that for  $k$ large, $\E(\Psi_{k})\leq 
c_{(\delta_{0},0)}/2$. This and \eqref{convforte} give the following uniformly a priori lower bound
$$
\|\psi_{k,1}\|_{2p}^{2p}>\sigma_{0},\qquad \sigma_{0}>0.
$$
This, \eqref{convforte} and \eqref{scalarc} allow us to argue as in \cite{cl}
(see also \cite{caz}), and by means of concentration compactness technique,
get the strong convergence (up to a subsequence) of $\psi_{k,1}$. Then,
there exists $\psi_{1}$ such that
$$
\|\psi_{1}\|_{2}^{2}=\|z^{\om}_{0}\|_{2}^{2},
$$
so that $\E_{0}(\psi_{1})\geq c_{\delta_{0}}$, but, passing to the limit 
in \eqref{scalarc} give  $\E_{0}(\psi_{1})= c_{\delta_{0}}$.
Therefore, there exist $\theta\in \R$ and $y\in \R^{n}$ such that
$\psi_{1}=e^{i\theta}z^{\om}_{0}(\cdot-y)$. Moreover, \eqref{consescal},
\eqref{gradient} and\eqref{equac} imply that $\psi_{k,2}\to0 $ in $H^{1}$, giving a contradiction with \eqref{absurd}.
\edim


{\bf Proof of Theorem \ref{mainstab3}.}\quad 
Arguing by contradiction, as in the proofs of the other stability results, 
and using Theorem \ref{carabound} and Proposition \ref{caraA}
we get a positive number $\eps_{0}$ and a sequence $\Psi_{k}$ 
such that
\beq\label{ultima}
\inf_{\theta\in\R,y\in \rn}\|\Psi_{k}
-(e^{i\theta_{1}}z^{\om}_{\bt}(\cdot-y),e^{i\theta_{2}}z^{\om}_{\bt}(\cdot-y))\|\geq \eps_{0},\quad \|\psi_{k,i}\|_{2}\to
\|z^{\om}_{\bt}\|_{2},\quad \E(\Psi_{k})\to c_{(\delta,\delta)}
\eeq
Following the proof of Lemma 2.3 in \cite{oh}, it can be proved that the 
sub-additivity condition holds for Problem \eqref{ohta}, then by 
concentration-compactness arguments (see Section IV in \cite{lions})
we obtain that $\Psi_{k}$ is compact and the conclusion follows passing to the limit in \eqref{ultima}.
\edim

\subsection{Instability in the critical or supercritical case}\label{subcrit}

In this subsection we will prove Theorem \ref{instap} and \ref{instapB}.

{\bf Proof of Theorem \ref{instap}.}\quad 
In order to prove conclusion ${\it a)}$ let us assume $1+2/n<p<n/(n-2)$
and consider first the set $\G$. Let $U\in \G$  so that 
$U\in {\mathcal P}$. When we fix $U_{s}=U^{s^{n/2},s}$ with $s>1$
we get that $\lambda_{*}(U_{s})<1$ and Conclusion $c)$ of Lemma 
\ref{glambda} implies that $\mathcal{R}(U_{s})<0$ and Conclusion
{\it b)} of Lemma \ref{glambda} gives 
\beq\label{questa}
\I(U_{s})=g(1)<g(\lambda_{*})=\I(U)=\mnehari.
\eeq
Let $\Phi_{s}$ the solution generated by $U_{s}$. 
By \eqref{energy} we have $\I(\Phi_{s})=\I(U_{s})$, when 
the solution exists. By continuity $\mathcal{R}(\Phi_{s}(t))<0$ 
for $t$ small; moreover, Lemma \ref{RI} and \eqref{questa} 
imply that
$$
\mathcal{R}(\Phi_{s}(t))\leq \I(\Phi_{s}(t))-\mnehari=-\sigma<0,
$$
showing that $\mathcal{R}(\Phi_{s}(t))<0$ when the solution exists. 
Defining the variance function
$$
V(t)=\||x|\Phi_{s}(t)\|_{2}^{2},
$$
it follows that $V^{''}(t)=8\mathcal{R}(\Phi_{s}(t))\leq -8\sigma$. 
Thus, there exists $T^{*}$ such that $V(T^{*})=0$ showing, by using
Hardy's  inequality, that $\Phi_{s}$ blows up in $T^{*}$ (see \cite{fm}), 
which gives conclusion {\it a)} for the set $\G$.
\\
Let now $U\in \Se$ then $U=(u_{1},0)$ (for example) with 
$u_{1}=e^{i\theta}z^{\om_{1}}_{0}$. $u_{1}$ is unstable for the single equation(with $\bt=0$) then Theorem 8.2.2 in \cite{caz} 
implies that there exists $u^{0}_{\eps,1}$ such that 
$\|u^{0}_{\eps,1}-u_{1}\|\leq \eps$ and the solution generated by $u^{0}_{\eps,1}$,  $\phi_{\eps,1}$ blows up in finite time. Now if we choose $U_{\eps}=(u^{0}_{\eps,1},0)$ we get $\|U_{\eps}-U\|\leq \eps $ and the solution generated by $U_{\eps}$ is (by the well posedness of the Cauchy problem) $\Phi_{\eps}=(\phi_{\eps,1},0)$ blows up in finite time.
\\
Now consider $p=1+2/n$. From Proposition \ref{infs} we get that any 
$U$ solution of \eqref{ellittico} satisfies \eqref{tre}, then
we get $\mathcal{R}(U)=0$ so that $\mathcal{R}(\la U)<0$ for
any $\la>1$. Let $U_{\la}=\la U$ be the initial datum of \eqref{schr}
and $\Phi_{\la}$ the corresponding solution. \eqref{energy} implies
\bdm
0>\mathcal{R}(U_{\la})=2\E(U_{\la})=
2\E(\Phi_{\la})=\mathcal{R}(\Phi_{\la}),
\edm
so that, also in this case, the variance is concave and
the solution $\Phi_{\la}$ blows up in finite time.
\edim
{\bf Proof of Theorem \ref{instapB}.}\quad 
Let $U\in \B$ then $U=(u_{1},u_{2})$ with $u_{1}=e^{i\theta_{1}}\zbo$ and $u_{2}=e^{i\theta_{2}}\zbo$, again $\zbo$ is unstable because it is a ground state for the single equation with coefficient $\bt+1$ in the nonlinear term, so that we can apply, as in the previous result, Theorem 8.2.2 in
\cite{caz} to obtain an initial datum $u_{\eps}$ such that
$\|u_{\eps}-\zbo\|\leq \eps$ and the solution that starts from $u_{\eps}$
$\phi_{\eps}$ blows up in finite time. If we choose $U_{\eps}=(u_{\eps},u_{\eps})$ we have (by the well posedness 
of the Cauchy problem) that the solution generating from $U^{\eps}$ 
is $\Phi_{\eps}=(\phi_{\eps},\phi_{\eps})$ and it blows up in finite
time. 
\edim

\section{Conclusions}

In summary, we have studied the problem of the orbital 
stability of standing waves in two weakly coupled 
nonlinear Schr\"odinger equations. 
In analogy of what happens for the single equation case we have that 
least action solutions give rise to orbitally stable standing waves. But, 
the system admits also other families of orbitally stable standing waves, for example the set of solutions with one trivial component whose elements are not
ground states (and have Morse index greater than one) for $\beta$ large.
Moreover, for $\om_{1}=\om_{2}$ least action solution with both nontrivial components also generate orbitally stable solutions of \eqref{schr}, and this  holds  also when they are not ground states of \eqref{ellittico}.
So that it seems that having Morse index, with respect to $\I$, equal to one is not a necessary property to gain orbital stability. In our opinion this is linked to the facts that the $L^{2}$ norms of the components are conserved separately.
\\
We remark that it remains open the question of the stability for the set of minima of the energy whose components have different $L^{2}$ norms,
at least for $\beta$ small. Moreover, it is an interesting open problem
to find conditions, maybe related to the geometrical properties of $\I$,
on a solution in order to produce instability. More precisely, it would be interesting to understand how to extend the result of \cite{gss} for this kind
of system.
\medskip

{\bf Acknowledgement.}\quad Part of this work was developed 
while the second autor was visiting the {\it Universidade de 
Bras\'\i{}lia}: he is very grateful to all the {\it Departmento 
de Matem\'atica} for the warm hospitality. Moreover he wishes 
to thank also professor Gustavo Gilardoni for many stimulating 
and friendly discussions.

The authors wish to thank Louis Jeanjean and Hichem Hajaiej  
for some very useful and stimulating observations.


\end{document}